\documentclass[11pt]{article}
\usepackage{amssymb,epsfig}
\usepackage{graphicx}
\usepackage{amsmath}
\usepackage{amsfonts}
\usepackage[latin1]{inputenc}

\def\R{{\Bbb R}}
\def\N{{\Bbb N}}

\usepackage{color}
\font\msbm=msbm10 scaled 1200

\def\R{\hbox{\msbm R}}

\newtheorem{defi}{\large\bf Definition}[section]
\newtheorem{theo}[defi]{\large\bf Theorem}
\newtheorem{lemma}[defi]{\large\bf Lemma}
\newtheorem{remark}[defi]{\large\bf Remark}

\title{ \bf
{On a Moser-Steffensen type method for nonlinear systems of equations}}
\author{
{ S. Amat$^a$, M. Grau-Sanchez$^b$, M. A. Hern\'{a}ndez-Ver\'{o}n$^c$, M. J. Rubio$^c$}\\
{\quad}
\\
{\small
$^a$ Dept. Applied Mathematics  and Statistics,}\\
{\small Polytechnique University of  Cartagena.
30230. Cartagena, Spain.}\\
{\small Research supported by the Fundaci\'on Seneca dentro del}\\
{\small Programa Jiménez de la Espada de Movilidad, Cooperaci\'on e Internacionalizaci\'on}\\
{\small $^b$ Dpt. of Applied Mathematics II,}\\
{\small Technical University of Catalonia.
08034 Barcelona, Spain.}\\
{\small Research supported by MTM2011-28636-C02-01 of the Spanish}\\
{\small Ministry of Science and Innovation. (Spain).}\\
$^c${\small
Dept. Mathematics and Computation,}\\
{\small University of La Rioja.
26004 Logro\~no, Spain.}\\
{\small Research supported by MTM2014-52016-C2-1-P of the Spanish}\\
{\small Ministry of Science and Innovation. (Spain).}\\
{\small E-mail:
 sergio.amat@upct.es, miquel.grau@upc.edu, mahernan@unirioja.es,}\\
{\small mjesus.rubio@unirioja.es
}}
\date{ }

\begin{document}

\maketitle
%%%%%%%%%%%%%%%%%%%%%%%%%%%%%%%%%%%%%%%%%%%%%%
\begin{abstract}
This paper is devoted to the construction and analysis of a Moser-Steffensen iterative scheme.
The method has quadratic convergence without evaluating any derivative nor inverse operator.
We present a complete study of the order of convergence for systems of equations, hypotheses ensuring
the local convergence and finally we focus our attention to its numerical behavior. The conclusion is that
the method improves the applicability of both Newton and Steffensen methods having the same
order of convergence.
\end{abstract}

\begin{flushleft}
 {\bf Classification A.M.S. :}   65J15, 47H17.\\
\end{flushleft}

\section{Introduction}

One of the most studied problems in numerical analysis is the solution of nonlinear equations
\begin{equation}\label{NewtEq}
F(x)=0.
\end{equation}

Many scientific and engineering problems can be written in the form of a nonlinear systems of equations, where $F$ is a nonlinear operator defined on a non-empty open convex subset $\Omega$ of $\R^m$ with values in $\R^m$. A powerful tool to solve these nonlinear systems of equations is by means of iterative methods.

It is well-known that Newton's method,

\begin{equation}\label{Newt}
   \left\{ \begin{array}{l} x_{0}\mbox{ given in} \, \, \Omega, \\[1ex]   x_{n+1}=x_n-F'(x_n)^{-1}F(x_n),\quad n\geqslant 0,
\end{array} \right.
\end{equation}
is the one of the most used iterative method to approximate the solution $x^*$ of (1). The quadratic convergence and the low operational cost of (2) ensure that Newton's method has a good computational efficiency. But the existence of the operator $F'(x_n)^{-1} $ is needed at each step or equivalently to solve $F'(x_n)(x_{n+1}-x_{n})=-F(x_{n})$. It is also known another feature of Newton's iteration: calculating inverse operators is not needed if the inverse of an operator is approximated; namely, to approximate $P^{ -1}$ , where $P \in \mathcal{L}( X , Y )$, the set of bounded linear operators from X into Y , we can phrase $G(Q) = Q^{-1} -P = 0$, so that Newton's method is reduced to
$Q_{n+1}=2Q_{n}-Q_{n}PQ_{n}, n\geqslant0$, provided that $Q_{0}$ is given. So, in \cite{Moser}, Moser proposed, to solve nonlinear operator equations, the following iterative process:

$$
\left\{ \begin{array}{l}
 x_{0}, \  B_0 \quad {\rm given},
\medskip
\\
x_{n+1}=  x_{n} - B_n F(x_{n}),\quad  n\geq 0
\medskip
\\
B_{n+1} = 2 B_{n} -  B_n F'(x_n) B_n, \quad  n\geq 0.
\end{array} \right.
$$

This new iterative method, which can be considered as a Newton-type method, does not need to calculate inverse
operators. It can be shown that the rate of convergence is  $\frac{1\sqrt{5}}{2}$, provided the root is simple (see \cite{Moser}).
This process uses the same amount of information per step as Newton's method, but it converges no faster than the secant method, so that this unsatisfactory from a numerical point of view. For that, Hald proposes in \cite{Hald} the use of the following iteration:

\begin{equation} \label{NewtUlm2}
\left\{ \begin{array}{l}
 x_{0}, \  B_0 \quad {\rm given},
\medskip
\\
x_{n+1}=  x_{n} - B_n F(x_{n}),\quad  n\geq 0
\medskip
\\
B_{n+1} = 2 B_{n} -  B_n F'(x_{n+1}) B_n, \quad n\geq 0.
\end{array} \right.
\end{equation}

Observe that the first equation is similar to the Newton's method in which case $B_n$ is $F'(x_n)^{-1}$. The second equation is Newton's method applied to $G(P) = P^{-1}- F'(x_{n+1}) = 0$.
We can stress two features of (4): it has the same rate of convergence of Newton's method and it does not need to solve nonlinear equations at each step (it is ``inversion free"). Moreover, from the convergence of the sequence $ {x_n}$, the convergence of the sequence ${B_n}$ can be given, so that it also produces successive approximations $B_n$ to the bounded right inverse of $F'(x^*)$ at the solution $x^*$, which is very helpful when sensitivity of solutions to small perturbations are investigated. However, as Newton's method, this method has the same serious shortcoming: the derivative $F'(x)$ has to be evaluated at each iteration. This makes it unapplicable to equations with nondifferentiable operators and in situations when evaluation of the derivative is too costly.
The goal of this paper is to modify the previous two-step method in order to avoid the evaluation of any Fr\'echet derivative.

It is common to approximate derivatives by divided differences, so that iterative methods that use divided differences instead of derivatives are obtained. Remember that an operator $[u,v;F]$, $u,v\in\Omega$, is called a first order divided difference  \cite{Potra,PoPt}, if it is a bounded linear operator such that
$$
[u,v;F]:\Omega\subset  \R^m  \longrightarrow  \R^m
\quad\text{and}\quad
[u,v;F](u-v) = F(u)-F(v).
$$

In this paper,  we will start with Steffensen's method \cite{AABL,Arg1,Arg2,ArgMag,Cord2,CHMT,Cord1,EHRV,GND,PID} to approximate a solution of the nonlinear system of equations  $F(x)=0$. It is a one-point iterative process given by the following algorithm:

\begin{equation}\label{metStef}
 \left\{ \begin{array}{l} x_{0}\mbox{ given}, \\[1ex] x_{n+1}=x_n - [x_{n},x_n +F(x_n);F]^{-1}F(x_n),  \, \, n\geq 0, \end{array} \right.
\end{equation}
The main interest  of this iterative process lies in that the approximation of the derivative $F'(x_n)$, that appears in each step of Newton's method, by the divided difference $[x_{n},x_n +F(x_n);F]^{-1}$ is good enough to keep the quadratic convergence \cite{Stef} and, therefore, keeps the computational efficiency of Newton's method.
However, in the Steffensen method remains an interesting problem to solve, its algorithm needs to solve a linear system of equations at each step:
\begin{equation}\label{SistStef}
 \left\{ \begin{array}{l} x_{-1}, x_{0}\mbox{ given}, \\[1ex]
[x_{n},x_n +F(x_n);F] (x_{n+1}-x_n)= F(x_n),  \, \, n\geq 0. \end{array} \right.
\end{equation}

In order to improve this fact, following the previous ideas applied in the case of Newton's method, we consider the following Moser-Steffensen type method

\begin{equation} \label{metUlm}
\left\{ \begin{array}{l}
 x_{0}, \  B_0 \quad {\rm given},
\medskip
\\
x_{n+1}=  x_{n} - B_n F(x_{n}),\quad  n\geq 0
\medskip
\\
B_{n+1} = 2 B_{n} -  B_n {\left[x_{n+1},x_{n+1} +F(x_{n+1});F\right]} B_n,
\end{array} \right.
\end{equation}
where we have changed the resolution of the linear system in every step by making several matrix multiplications.

Whereas the operational cost of both processes is similar, however the possibility of {\em ill-conditioned} linear systems appearing  will be avoided simply taking appropriate initial matrix $B_0$ , and therefore this previous algorithm will be able to be stable more easily. The matrix given by the divided difference $[x_{n},x_n +F(x_n);F]$, can be {\em ill-conditioned} and therefore the linear system of equations previously indicated (\ref{SistStef}), in the classical Steffensen method (\ref{metStef}),  should cause instabilities.

\begin{remark}
\medskip
Along the paper, we consider  the following first order divided difference in $ \R^m $:
$[u, v;F]=([u,v;F]_{ij})_{i,j=1}^{m}$,
where
{\small
\begin{equation}\label{DD}
    [u, v;F]_{ij} = \frac{1}{u_{j}-v_{j}}
\left(F_{i}(u_{1},\ldots,u_{j},v_{j+1},\ldots,v_{m})
- F_{i}(u_{1},\ldots,u_{j-1},v_{j},\ldots,v_{m})\right),
\end{equation}
}
$ u= (u_1,u_2,\dots,u_m)^T$ and $v=(v_1,v_2,\dots,v_m)^T.$
\end{remark}

The organization of the paper is the following. In Section 2 we study the local error equation, seeing his quadratic convergence. After, in Section 3, we obtain a local convergence result for this iterative process. Moreover, in the Section 4, the numerical behavior is analyzed. From this study, we can conclude that the considered method (\ref{metUlm}) improves the applicability of Steffensen method (\ref{metStef}).

%%%%%%%%%%%%%ERROR LOCAL%%%%%%%%%%%%%%%%%%%%%%%%%%
\section{Local order of convergence}

We assume that $ F:\Omega\subseteq {\mathbf{R}^m} \longrightarrow {\mathbf{R}^m}$ has at least $2$-order derivatives with continuity on $\Omega$ for any $x\in {\Omega}$  lying in a neighborhood of a simple zero, $x^*\in \Omega$, of the
system $\, F (x) = 0$. We can apply Taylor's formulae to $F(x)$. By setting $\,e_k = x_k - x^*$, the local order,
and assuming that $\, {\left[F^{\prime}\left(x^*\right)\right]}^{-1}$ exists, we have
\begin{equation}\label{fun01}
   F(x_k)\,=\, F(x^* + e_k)\,=\, \Gamma \left( \, e_k +\,O(e^2_k)\,\right),
\end{equation}
where $ \Gamma = F^{\prime}\left(x^*\right),  \;e_k^{2}\,=\,(e_k, e_k)\in \mathbf{R}^m \times \mathbf{R}^m$. Moreover, from  \cite{ggn1} we obtain
\begin{equation}\label{divdif01}
\left[x_k\,, y_k\,;\;F\right] =  \Gamma \left( I + A_2(e_k+\varepsilon_k) + O_2(e_k,\varepsilon_k) \right) ,
\end{equation}
where $\varepsilon_k=y_k-x^*$ and $A_2 = \frac{1}{2}\,\Gamma^{-1}\,F^{\prime \prime}(x^*)\in \mathfrak{L}_2(\Omega,\mathbf{R}^m)$.
We say that an operator depending on $e_k$ and $\varepsilon_k$ is an $O_{2}(e_k, \varepsilon_k)$ if it is an
$O(e^{q_0}_k\,\varepsilon^{q_1}_k)$ with $q_0+q_1=2\,,$ $q_i\ge 0\,,\ i=0,1$.

\vspace{1mm}
If we expand in formal power series of $e_k$ and $\epsilon_k$, the inverse of the divided difference given in \eqref{divdif01} is:

\vspace{-2mm}
\begin{equation} \label{divdif02}
\left[x_k\,,y_k\,;\,F\right]^{-1} = \left( I -A_2( e_k + \epsilon_k) +  O_{2}(e_k, \varepsilon_k)\right) {\Gamma}^{-1}.
\end{equation}
These developments of the divided difference operator \eqref{divdif02} were previously used in
the study of Grau~{\sl et al.}~\cite{ggn1}. In particular, for the operator $\Theta_k\,=\, \left[x_k\,,x_k + F(x_k)\,;\,F\right]$ we take
$$
\left[x_k\,,x_k + F(x_k)\,;\;F\right] =  \Gamma \left( I + A_2(e_k+ ( 2\,I+\Gamma) \,e_k) + O(e^2_k) \right) ,
$$
Note that we have the following relation between $F(x_k)$ and $\Theta_k\,e_k$:
\begin{equation}\label{r0}
 F(x_k) \,=\, \Theta_k\,e_k + O(e^2_k) .
\end{equation}

A first main local result is stated in the following theorem:

\vspace{3mm}
\begin{theo}\quad $e_{k+1}\,=\,E_k \,e_k + O(e^2_k)$, \quad {\rm where} $\; \: \, E_k \,= \, I - B_k\, \Theta_k$.
\end{theo}

\textit{\textbf{Proof.}}\quad

 By subtracting $x^*$ from both sides of the first step of
 (\ref{metUlm}) and  taking into account (\ref{r0}) we obtain
\begin{eqnarray} \nonumber
e_{k+1} &=& e_k - B_k\, F(x_k) \\ \nonumber
        &=& e_k - B_k\,\Theta_k \,e_k + O(e^2_k) \\ \label{prop1}
        &=& E_k \,e_k + O(e^2_k) ,
\end{eqnarray}
and the proof is completed.
\qquad$\blacksquare$
\bigskip

%\end{proof}%%%%%%%%%%%%%%%%%%%%%%%%%%%%%%%%%%%%%%%%%%%%%

\vspace{3mm}
Other relations that are necessaries to obtain a second theorem are the following ones:
\begin{eqnarray} \nonumber
E_{k+1} &=& I - B_{k+1}\, \Theta_{k+1} \\ \nonumber
        &=& I - 2\,B_k\,\Theta_{k+1} + \,B_k\,\Theta_{k+1}\,B_k\,\Theta_{k+1} \\ \label{r1}
        &=& \left( I - \,B_k\,\Theta_{k+1} \right)^2.
\end{eqnarray}
and
\begin{eqnarray} \nonumber
I - \,B_k\,\Theta_{k+1} &=& I - B_k\,(\Theta_{k+1} -\Theta_k)\, B_k\,\Theta_{k} \\ \nonumber
        &=& E_k - B_k\,(\Theta_{k+1} -\Theta_k) \\ \label{r2}
        &=& E_k - G_k\,e_k + O(e_{k+1}, e^2_k) ,
\end{eqnarray}
where $\, G_k = B_k\,A_2\,(2I+\Gamma) \in \mathfrak{L}_2 (\Omega, \mathbf{R})$.

\vspace{2mm}
A second main local result is stated in the following theorem:

\vspace{1mm}
\begin{theo}\quad
$E_{k+1}\,=\,\left(E_k + G_k\,e_k\right)^2 + O_2(E_k, e_k)$.
\end{theo}

\textit{\textbf{Proof.}}\quad

 By substituting (\ref{r2}) into (\ref{r1}) we have
\begin{eqnarray} \nonumber
E_{k+1} &=& \left( E_k - G_k\, e_k) + O( E_k\,e_{k}, e^2_k) \right)^2 \\ \label{prop2}
        &=& E^2_k + \left( G_k\,e_k\right)^2 - E_k\,G_k\,e_k - G_k\,e_k\, E_k + o(E^2_k, E_k\,e_k) .
\end{eqnarray}
The proof is complete.
\qquad$\blacksquare$
\bigskip

%\end{proof}%%%%%%%%%%%%%%%%%%%%%%%%%%%%%%%%%%%%%%%%%%%%%

\vspace{3mm}
In a more precise way we can write Theorems 1 and 2 (see (\ref{prop1}) and (\ref{prop2}) respectively) in norm terms. Namely,
\begin{eqnarray} \label{A1}
\|e_{k+1}\| &=& O \left( \| E_k \| \|e_k\| \right) \\[0.5em] \label{A2}
\|E_{k+1}\| &=& O_2\left( \| E_k \| ,\|e_k\| \right) = O\left( \| E_k \|^2 ,\|e_k\|\|^2, \| E_k \| \|e_k\|\right) .
\end{eqnarray}

We have three possibilities:
\begin{itemize}
\item If  $\|E_{k+1}\| = O \left( \| E_k \|^2 \right)$, then  $\|E_{k}\| = O \left( \| E_{k-1} \|^2 \right)$, and from (\ref{A1}) we take $\|e_{k+1}\| = O \left( \| E_{k-1} \|^2 \|e_k\| \right)$, and applying (\ref{A1}) for $k-1$, we have $\|e_{k}\| = O \left( \| E_{k-1} \| \|e_{k-1}\| \right)$. Hence, $\|e_{k+1}\| = O \left( \| E_{k-1} \|^3 \|e_{k-1}\| \right)$. Taking into account that from (\ref{A1}) with $k-1$ we have $\|e_{k}\|^3 = O \left( \| E_{k-1} \|^3 \|e_{k-1}\|^3 \right)$, finally, we obtain
    \begin{equation}\label{edd1}
         \|e_{k+1}\| = O \left( \| e_{k} \|^3 \|e_{k-1}\|^{-2} \right) .
    \end{equation}
The  equation associated with (\ref{edd1}) is
$\, p_1(t)=t^{2} - 3 t +2 = 0$. The  real positive root that coincides with the local order of
convergence in this case is  $\rho = 2$  (see \cite{OrRh,Torn}).

\item If  $\|E_{k+1}\| = O \left( \|e_k\|\|^2 \right)$ and $E_{k} = O \left( \|e_{k-1}\|\|^2 \right)$, then from (\ref{A1}) we have
 \begin{equation}\label{edd2}
         \|e_{k+1}\| = O \left( \| e_{k} \| \|e_{k-1}\|^{2} \right) .
    \end{equation}
In this case, the local order of convergence $\rho$ is the unique real positive root of the indicial polynomial (see \cite{OrRh,Torn,Trau}) of the error difference equation \eqref{edd2} given by $p(t)=t^2-t-2$. That is, $\rho=2$.

\item If  $\|E_{k+1}\| = O \left( \| E_k \| \|e_k\|\right)$ and we have (see (\ref{A1})) $\|e_{k+1}\| = O \left( \| E_k \| \|e_k\| \right)$, then
\begin{eqnarray*}
\|e_{k+1}\|\;=\; O \left( \|E_{k+1}\| \right) &=& O \left( \| E_k \| \| E_{k-1} \|\|e_{k-1}\|\right) \\[0.5em]
 &=& O \left(\|E_{k-1}\|^2 \|e_{k-1}\|^2 \right) \\[0.5em]
 &=& O \left( \|e_{k}\|^2 \right) .
\end{eqnarray*}
The method, as in the two first cases, presents quadratic order of convergence again.
\end{itemize}

We are now in a position to state the following theorem.

%\vspace{1mm}
 \begin{theo}
 \quad
The iterative method defined by (\ref{metUlm}), from a local view of point, is a quadratic method. That is $ \|e_{k+1}\|\;=\; O \left( \|e_{k}\|^2 \right) $. \quad\quad\quad\quad$\blacksquare$
\end{theo}

\vspace{5mm}

%%%%%%%%%%%%%%%%%%%%%%%%%%%%%%%%%%%%%%%%%%%%%%%%%%%%%%%%%%%
\section{Local convergence analysis}
\label{SC}
%%%%%%%%%%%%%%%%%%%%%%%%%%%%%%%%%%%%%%%%%%%%%%%%%%%%%%%%%%%

In this section, we prove the local convergence of the Moser-Steffensen type method given by (\ref{metUlm}).
First, we establish a system of recurrence relations, from the real parameters that are introduced under some conditions for the pair $(F,x_{0})$, where a sequence of positive real numbers is involved.
After that, we can guarantee the semilocal convergence of the method in $\R^m$.

%%%%%%%%%%%%%%%%%%%%%%%%%%%%%%%%%%%%%%%%%%%%%%%%%%%%%%%%%%%
\subsection{Recurrence relations}
%%%%%%%%%%%%%%%%%%%%%%%%%%%%%%%%%%%%%%%%%%%%%%%%%%%%%%%%%%%

We suppose:
%\begin{quote}

\begin{itemize}

\item [(C1)]%
Let $x^*$ a zero of nonlinear system $F(x)=0$ such that $ \Vert F'(x^*)\Vert \le M$ and there exists $r>0$ with $B(x^*,r)\subseteq \Omega$.

\item [(C2)]%
Let  $ B_0\in{\cal L}(\R^m,\R^m)$ with ${B_0\ne 0}$  such that  $\Vert { B_0}\Vert = \beta$ and $\Vert I-{ B_0\, F'(x^*)}\Vert = \delta< 1$.

\item [(C3)]%
Consider $x_0\in B(x^*,r)$, such that there exists $\tilde{r}>0$ such that $x_0+F(x_0) \in B(x^*, \tilde{r})$, with $\in B(x^*, \tilde{r}) \subseteq \Omega$.

\item [(C4)]%
 For each pair of distinct points $x, y \in \Omega $, there exists a first order
divided difference ${\left[ x,y; F \right]}$ of $F$ and  $k \geqslant  0$, such that
\begin{equation}\label{centr-divdif}
\Vert { {\left[ x, y; G \right]}- F'(x^*)}\Vert
\le
k {\left(  {\Vert x-x^*\Vert} +{\Vert y-x^*\Vert}\right)}; \ \
\forall x, y \in \Omega,
\end{equation}

\end{itemize}

%\end{quote}

Notice that, in these conditions, the Fr\'echet derivative of $F$ exists in $\Omega$ and satisfies
$ {{\left[ x, x; F \right]} = F'(x)}$. On the other hand, to simplify the notation,  we denote $L_n=\displaystyle{   {\left[ x_{n}, \, x^*; F\right]}  }$ and ${ \Theta_n} = { \left[ x_{n}, x_{n}+F(x_n); F\right] }$, for all $n \in \N$.
\medskip

From the above, we denote $\alpha_0=r$, $\beta_0=\beta$, $\delta_0=\delta$, $\tilde{\alpha}_0= \tilde{r}$ and define the scalar sequences:
\begin{equation}\label{realsucs}
 \begin{array}{l}
 \alpha_n= (\delta_{n-1} + k \,  \beta_{n-1} \alpha_{n-1}) \alpha_{n-1}, \quad
   \tilde{\alpha}_{n} = (1 + M +  k \,   \alpha_{n}) \alpha_{n}
    \medskip
    \\
d_n= \delta_{n} + k \,  \beta_{n}  (\alpha_{n+1} +  \tilde{\alpha}_{n+1}),  \quad  \beta_n=(1+d_{n-1})\beta_{n-1}
   \medskip
\\
\delta_n= \delta^2_{n-1} +  k \,  M \beta^2_{n-1}  (\alpha_{n} +  \tilde{\alpha}_{n})
\end{array}
\end{equation}

Next, for $n=1$, we prove he following recurrence relations for sequences (\ref{realsucs}) and $\{x_{n}\}$:
\begin{align}  \nonumber
\Vert { x_1} - x^*\Vert  & < \alpha_1, \\ \nonumber
\Vert { x_1}+F(x_1) - x^*\Vert  &< \tilde{\alpha}_{1} ,\\ \nonumber
\Vert  I -   B_1 { \Theta_{2}}\Vert  &\le d_1, \\ \nonumber
 \Vert B_{1}  \Vert   &     \le \beta_1, \\ \nonumber
  {\Vert I - B_1\,  F'( x^*)\Vert}  &  \le \delta_1, \nonumber
\end{align}
provided that
$$
x_0, x_0+F(x_0)\in\Omega.
$$

From the initial hypotheses, it follows that $x_1=  x_0 -B_0 F(x_{0})$ is well defined and
\begin{align}\nonumber
\Vert { x_1} - x^*\Vert
&= \Vert x_0-{ B_0} F(x_{0}) - x^*\Vert \\ \nonumber
& = \Vert x_0- x^*- { B_0}  L_0  (x_0- x^*) \Vert   \\ \nonumber
& \le \Vert I - { B_0} {L_0} \Vert  \Vert x_0- x^* \Vert  \\ \nonumber
& \le \left( {\Vert I - B_0\,  F'( x^*)\Vert}  +  {\Vert B_0 \Vert}   {\Vert L_0 -  F'( x^*)\Vert} \right) {\Vert x_{0} - x^*\Vert} \\ \nonumber
& < \left(\delta_0  +  \beta_0 \,  k  \alpha_0  \right) {\Vert x_0 -  x^*\Vert} \\ \nonumber
& < \left(\delta_0  +  \beta_0 \,  k  \alpha_0\right)   { \alpha_0}= { \alpha_1}.\\ \nonumber
 \end{align}

 On the other hand,
\begin{align}\nonumber
\Vert   x_{1}+F(x_1) - x^* \Vert
& =  {\Vert x_{1}+ L_1\,  (x_1 - x^* ) -x^*\Vert}  \le   {\Vert I+ L_1\Vert } {\Vert x_1 - x^*\Vert} \\ \nonumber
&  \le \left({\Vert I+F'( x^*)\Vert}  + {\Vert F'( x^*)- L_1 \Vert }  \right) {\Vert x_{1} - x^*\Vert} \\ \nonumber
&  \le \left( {1+\Vert F'( x^*)\Vert}  + k\, {\Vert x_{1} - x^*\Vert } \right) {\Vert x_{1} - x^*\Vert} \\ \nonumber
& \le  \left( {1+M}  + k\, {\Vert x_{1} - x^*\Vert } \right) {\Vert x_{1} - x^*\Vert} \\ \nonumber
&  \le  \left({1+M}  + k\,  \alpha_{1}\right) \alpha_{1} = \tilde{\alpha}_{1}. \\ \nonumber
 \end{align}

 Assuming  that ${\alpha_1} < {\alpha_0}$ and $(1+M+k r) r< \tilde{r}$, then $ \tilde{\alpha}_{1}< \tilde{\alpha}_{0}$ and therefore
$ x_{1}, x_{1}+F(x_1)\in {\Omega }$. So,  there exist $ { \Theta_1} = { \left[ x_{1}, x_{1}+F(x_1); F\right] }$
and   $B_1=   2 B_{0} -  B_0 { \Theta_1} B_0$.
Then, we establish
\begin{align}\nonumber
\Vert  I -   B_0  { \Theta_1}\Vert
& \le \Vert  I -   B_0  {F'(x^*)}\Vert  + \Vert     B_0 \Vert  \Vert   { \Theta_1} - {F'(x^*)} \Vert  \\ \nonumber
& \le \Vert  I -   B_0  {F'(x^*)}\Vert  + \Vert     B_0 \Vert  \Vert \,  k (   {\Vert x_{1} - x^*\Vert} +  {\Vert x_{1} +F(x_1)- x^*\Vert}  )  \\ \nonumber
& \le \delta_0  + \beta_0 \, k (   {\alpha_1} +  \tilde{\alpha}_{1} ) = d_0. \\ \nonumber
 \end{align}
As a consequence,
\begin{align}\nonumber
  \Vert B_{1}  \Vert
 & = \Vert 2 B_{0} -  B_0  { \Theta_1} B_0 \Vert \\ \nonumber
 &   \le  \left(  1+ \Vert  I -   B_0  { \Theta_1}\Vert   \right)  \Vert  B_0  \Vert  \\ \nonumber
&   \le \left(  1+ d_0   \right) \Vert  B_0  \Vert = \beta_1. \\ \nonumber
 \end{align}
Moreover, to finish the first step, notice that
\begin{align}\nonumber
{\Vert I - B_1\,  F'( x^*)\Vert}
& =  {\Vert I - \left(2 B_{0} - B_{0} \Theta_{1}  B_{0}\right) \,  F'( x^*)\Vert}  \\ \nonumber
& \le {\Vert I - B_{0}\,  F'( x^*)\Vert}^2 +  {\Vert B_{0} \Vert}\,  {\Vert \Theta_{1}-  F'( x^*)\Vert}  \,  {\Vert B_{0}\, F'( x^*)\Vert}  \\ \nonumber
& \le  {\Vert I - B_{0}\,  F'( x^*)\Vert}^2 +  {\Vert B_{0} \Vert}^2\, {\Vert F'( x^*)\Vert} \,
 k\, \left( {\Vert x_{1} - x^* \Vert  + \Vert x_{1}+ F( x_1)-x^*\Vert} \right)  \\ \nonumber
&< {\delta}^2_0 +  {\beta}^2_0 \, M \, k\, \left( \alpha_1  + \tilde{\alpha}_{1}\right) = {\delta_1}.  \\ \nonumber
 \end{align}

Next, for $n=2$, we prove the following recurrence relations for sequences (\ref{realsucs}) and $\{x_{n}\}$:

\begin{align}\nonumber
\Vert { x_2} - x^*\Vert  & < \alpha_2, \\ \nonumber
\Vert { x_2}+F(x_2) - x^*\Vert  &< \tilde{\alpha}_{2} , \\ \nonumber
\Vert  I -   B_2  { \Theta_{3}}\Vert  &\le d_2, \\ \nonumber
 \Vert B_{2}  \Vert   &     \le \beta_2, \\ \nonumber
  {\Vert I - B_2\,  F'( x^*)\Vert}  &  \le \delta_2.\nonumber
\end{align}

Now, by (\ref{metUlm}), $x_{2}=  x_{1} - B_1 F(x_{1})$, furthermore
\begin{align}\nonumber
\Vert x_{2} - x^*\Vert
& \le {\Vert I - B_1\,  [x_1,  x^*;\, F] \Vert}  {\Vert x_{1} - x^*\Vert} \\ \nonumber
& = {\Vert I - B_1\,  L_1\Vert}  {\Vert x_{1} - x^*\Vert}  \\ \nonumber
& \le \left({\Vert I - B_1\,  F'( x^*)\Vert}  +  {\Vert B_1 \Vert}   {\Vert L_1 -  F'( x^*)\Vert}  \right) {\Vert x_{1} - x^*\Vert}. \\ \nonumber
& \le \left(\delta_1 +  \beta_1 \,  k \, \alpha_{1} \right) \alpha_{1}  = {\alpha_2}. \\ \nonumber
 \end{align}
Moreover, it is easy to check
\begin{align} \nonumber
 {\Vert x_{2}+ F( x_2)-x^*\Vert}
& = {\Vert x_{2}+ L_2\,  (x_2 - x^* ) -x^*\Vert}  \le   {\Vert I+ L_2 \Vert } {\Vert x_2 - x^*\Vert} \\ \nonumber
& \le \left( {\Vert I+ F'( x^*)\Vert}  + {\Vert F'( x^*)- L_2 \Vert }  \right) {\Vert x_{2} - x^*\Vert} \\ \nonumber
&   \le \left({1+\Vert F'( x^*)\Vert}  + k\, {\Vert x_{2} - x^*\Vert } \right) {\Vert x_{2} - x^*\Vert} \\ \nonumber
& \le  \left(  {1+M}  + k\, {\Vert x_{2} - x^*\Vert }  \right) {\Vert x_{2} - x^*\Vert} \\ \nonumber
&   \le  \left({1+M}  + k\,  \alpha_{2} \right) \alpha_{2} = \tilde{\alpha}_{2}. \\ \nonumber
 \end{align}

Assuming  that $\alpha_{2} < \alpha_{1}$,  then $ \tilde{\alpha}_{2} < \tilde{\alpha}_{1} $ and therefore $x_2,x_2+F(x_2)\in \Omega$. So,
${ \Theta_2} $ and $B_2 $ are well defined.

Then, we have
\begin{align} \nonumber
\Vert  I -   B_1  { \Theta_2}\Vert
& \le \Vert  I -   B_1 {F'(x^*)}\Vert  + \Vert     B_1 \Vert  \Vert   { \Theta_2} - {F'(x^*)} \Vert \\ \nonumber
& \le \Vert  I -   B_1  {F'(x^*)}\Vert  + \Vert     B_1 \Vert  \Vert \, k (   {\Vert x_{2} - x^*\Vert} +  {\Vert x_{2} +F(x_2)- x^*\Vert}  ) \\ \nonumber
& \le \delta_1  + \beta_1\, k (   {\alpha_2} +  \tilde{\alpha}_{2} ) =  d_1, \\ \nonumber
\end{align}
and we get
\begin{align} \nonumber
  \Vert B_{2}  \Vert
& = \Vert 2 B_{1} -  B_1  { \Theta_2} B_1 \Vert  \\ \nonumber
 & \le  \left(  1+ \Vert  I -   B_1  { \Theta_2}\Vert   \right)  \Vert  B_1  \Vert  \\ \nonumber
 &  \le \left(  1+ d_1   \right) \Vert  B_1  \Vert \\ \nonumber
 & \le  \left(  1+ d_1   \right) \beta_1= \beta_2. \\ \nonumber
\end{align}

To finish the second step, we consider
\begin{align}\nonumber
{\Vert I - B_2\,  F'( x^*)\Vert}
& =  {\Vert I - \left(2 B_{1} - B_{1} \Theta_{2}  B_{1}\right) \,  F'( x^*)\Vert}  \\ \nonumber
& \le {\Vert I - B_{1}\,  F'( x^*)\Vert}^2 +  {\Vert B_{1} \Vert}\,  {\Vert \Theta_{2}-  F'( x^*)\Vert}  \,  {\Vert B_{1}\, F'( x^*)\Vert}  \\ \nonumber
& \le  {\Vert I - B_{1}\,  F'( x^*)\Vert}^2 +  {\Vert B_{1} \Vert}^2\, {\Vert F'( x^*)\Vert} \,
 k\, \left( {\Vert x_{2} - x^* \Vert  + \Vert x_{2}+ F( x_2)-x^*\Vert} \right)  \\ \nonumber
&< {\delta_1}^2 +  {\beta_1}^2\, M \, k\, \left( \alpha_2  + \tilde{\alpha}_{2}\right) = {\delta_2}.  \\ \nonumber
 \end{align}

At this time, we are able to obtain a general result that allows us to relate the sequences (\ref{realsucs}) and $\{x_{n}\}$.

\begin{lemma}
In the previous conditions, if the sequence $\{\alpha_n\}$ is decreasing and $(1+M+k r) r< \tilde{r}$, then
\begin{enumerate}
  \item[$(I)$] $\Vert { x_n} - x^*\Vert   < \alpha_n$,
  \item[$(II)$] $\Vert { x_n}+F(x_n) - x^*\Vert  < \tilde{\alpha}_{n}$,
  \item[$(III)$] $\Vert  I -   B_n  { \Theta_{n+1}}\Vert  \le d_n$,
  \item[$(IV)$] $ \Vert B_{n}  \Vert  \le \beta_n$,
   \item[$(V)$] $  {\Vert I - B_n\,  F'( x^*)\Vert}  \le \delta_n$,
\end{enumerate}
for all $n \in N$.
\end{lemma}

\textit{\textbf{Proof.}}\quad

First, since $\{\alpha_n\}$  is a decreasing sequence and $(1+M+k r) r< \tilde{r}$, then $\{\tilde{\alpha}_n\}$  is also decreasing. Therefore, we have that  $ x_{k}, x_{k}+F(x_k)\in {\Omega }$  is verified for $k=0,\ 1,\dots,\, n$, and therefore there exist ${ \Theta_k}, \, B_k$ such that ${x_{k+1}= x_k- B_k F(x_k)}$ is well defined.

Once shown the relationships $(I)-(V)$ in their first two steps previously, we will use mathematical induction. Suppose that the relations $(I)-(V)$ are true for  $k=1,\dots,\, n$ and we are going to show them for $k = n + 1$.

Observe that
$$
F(x_n) = F(x_n) -F(x^*)
= L_n\,  (x_n - x^* )
$$
and by (\ref{metUlm})
$$
x_{n+1} - x^* =  x_n -B_nF(x_n) -x^*
=  (I - B_n \,  L_n )(x_n -x^* ),
$$
 furthermore
$$
\Vert x_{n+1} - x^*\Vert
\le
{\Vert I - B_n\,  L_n\Vert}  {\Vert x_{n} - x^*\Vert}
$$
$$
\le \left(
 {\Vert I - B_n\,  F'( x^*)\Vert}  +  {\Vert B_n \Vert}   {\Vert L_n -  F'( x^*)\Vert}
\right) {\Vert x_{n} - x^*\Vert}
$$
$$
\le \left(
 {\Vert I - B_n\,  F'( x^*)\Vert}  +  {\Vert B_n \Vert}\,  k \, {\Vert x_{n} - x^*\Vert}
\right) {\Vert x_{n} - x^*\Vert}
\le \left(\delta_n +  \beta_n \,  k \, \alpha_{n}
 \right) {\Vert x_{n} - x^*\Vert}
$$
$$
\le \left(\delta_n +  \beta_n \,  k \, \alpha_{n}
 \right) {\alpha_{n}}=\alpha_{n+1}.
$$
In addition
\begin{align}\nonumber
 {\Vert x_{n+1}+ F( x_{n+1})-x^*\Vert} \nonumber
 & \le   {\Vert I+ L_{n+1} \Vert } {\Vert x_{n+1} - x^*\Vert}  \\ \nonumber
 & \le \left( {\Vert I+ F'( x^*)\Vert}  + {\Vert F'( x^*)- L_{n+1} \Vert } \right) {\Vert x_{n+1} - x^*\Vert}  \\ \nonumber
 & \le \left({1+\Vert F'( x^*)\Vert}  + k\, {\Vert x_{n+1} - x^*\Vert }  \right) {\Vert x_{n+1} - x^*\Vert} \\ \nonumber
 & \le  \left({1+M}  + k\, {\Vert x_{n+1} - x^*\Vert } \right) {\Vert x_{n+1} - x^*\Vert} \\ \nonumber
 &   \le  \left({1+M}  + k\, \alpha_{n+1}  \right) {\alpha_{n+1}} \\ \nonumber
 & = \tilde{\alpha}_{n+1}.  \nonumber
\end{align}

Assuming  that ${\alpha_{n+1}} < {\alpha_n}$, then  $ \tilde{\alpha}_{n+1}<\tilde{\alpha_n}$  and
$ x_{n+1}, x_{n+1}+F(x_{n+1})\in {\Omega }$. Consequently $\Theta_{n+1}$ and $B_{n+1}$ are well defined.

So, we can consider
\begin{align}\nonumber
\Vert  I -   B_n  { \Theta_{n+1}}\Vert  \nonumber
& \le \Vert  I -   B_n {F'(x^*)}  + \Vert     B_n \Vert  \Vert   { \Theta_{n+1}} - {F'(x^*)} \Vert \\ \nonumber
& \le \Vert  I -   B_n  {F'(x^*)}\Vert  + \Vert     B_n \Vert  \Vert \, k (   {\Vert x_{n+1} - x^*\Vert} +  {\Vert x_{n+1} +F(x_{n+1})- x^*\Vert}  ) \\ \nonumber
& \le \delta_n  + \beta_n\, k (   {\alpha_{n+1}} +  \tilde{\alpha}_{n+1} ) = d_n, \nonumber
\end{align}
that implies
$$
  \Vert B_{n+1}  \Vert
  = \Vert 2 B_{n} -  B_n  { \Theta_{n+1}} B_n \Vert
    \le  \left(  1+ \Vert  I -   B_n  { \Theta_{n+1}}\Vert   \right)  \Vert  B_n  \Vert
     \le
    \left(  1+ d_n   \right) \beta_n = \beta_{n+1}.
$$
Finally, to round off
$$
 {\Vert I - B_{n+1}\,  F'( x^*)\Vert} =  {\Vert I - \left(2 B_{n} - B_{n} \Theta_{n+1}  B_{n}\right) \,  F'( x^*)\Vert}
 $$
 $$
\le
 {\Vert I - B_{n}\,  F'( x^*)\Vert}^2 +  {\Vert B_{n} \Vert}\,  {\Vert \Theta_{n+1}-  F'( x^*)\Vert}  \,
 {\Vert B_{n}\, F'( x^*)\Vert}
$$
 $$
\le
 {\Vert I - B_{n}\,  F'( x^*)\Vert}^2 +  {\Vert B_{n} \Vert}^2\, {\Vert F'( x^*)\Vert} \,
 k\, \left(
  {\Vert x_{n+1} - x^* \Vert  + \Vert x_{n+1}+ F( x_{n+1})-x^*\Vert}
\right)
$$
$$
\le
 {\delta_{n}}^2 +  {\beta}^2_{n} \,M \,
 k\, \left(
  \alpha_{n+1}  + \tilde{\alpha}_{n+1}
\right)  =\delta_{n+1}.
 $$
The proof is complete.
\qquad$\blacksquare$
\bigskip

Once generalized the previous recurrence relations to every point of the sequence $\{x_n\}$, we have to guarantee that $\{x_n\}$ is a Cauchy sequence having into account these recurrence relations. For this, we first analyse the scalar sequences given by(\ref{realsucs}) in the next section.

%%%%%%%%%%%%%%%%%%%%%%%%%%%%%%%%%%%%%%%%%%%%%%%%%%%%%%%%%%%
\subsection{Analysis of the scalar sequence}
\label{se:ssWH}
%%%%%%%%%%%%%%%%%%%%%%%%%%%%%%%%%%%%%%%%%%%%%%%%%%%%%%%%%%%

Now, we analyse the scalar sequence defined in (\ref{realsucs}) in order to prove later the semilocal convergence of the sequence $\{x_{n}\}$ in  $\R^m$.
For this, it suffices to see that $\{x_{n}\}$ is a Cauchy sequence. First, we give a technical lemma.

 \begin{lemma}\label{Lem2}
 \quad Let $\{ \alpha_n\}$, $\{  \tilde{\alpha}_{n}\}$, $\{ d_n\}$ and $\{ \delta_n\}$
 be the  sequences  given by
(\ref{realsucs}).
If it is verified that
\begin{equation} \label{cond}
(1+M+k r) r< \tilde{r},
 \quad \delta_1 < \delta_0
 \quad {\rm and } \
\quad (1+d_0)^2 ( \delta_0 + k \beta \alpha_0 )<1,
\end{equation}
then
\begin{description}
\item[{(a)}]  \quad $  ( \delta_0 + k \beta_0 \alpha_0)<1$ and $  (1+d_0) ( \delta_0 + k \beta_0 \alpha_0 )<1$,
\item[{(b)}]  \quad the  sequences   $\{ \alpha_n\}$, $\{ {n}\}$, $\{ d_n\}$ and $\{ \delta_n\}$
are decreasing.
\end {description}
\end {lemma}
\textit{\textbf{Proof.}}\quad

 Observe  that as  $ (1 + d_0)> 1$, then
 ${  {(1+d_0)}^2 ( \delta + k \beta \alpha ) < 1}$ implies that
 $ (1+d_0) ( \delta + k \beta \alpha )<1$. By the same reason  $( \delta + k \beta \alpha )<1$ and (a)  holds.

We shall prove {\bf(b}) by induction.

From  {\bf(a)}, for $n=1$, we have that
$\alpha_1 = \left(\delta_0  + k\, \beta_0\,    \alpha_0\right) \alpha_0< \alpha_0$  and
$\tilde{\alpha}_{1}  =  \left({1+M}  + k\,  \alpha_{1} \right) \alpha_{1} <  \left({1+M}  + k\,  \alpha_{0} \right) \alpha_{0} =\tilde{\alpha}_{0}$.

For $n=2$, having into account that
$\beta_1  \alpha_1 <  (1+d_0)  \beta_0  \left(\delta_0  + k\, \beta_0\,    \alpha_0\right) \alpha_0  < \beta_0  \alpha_0 $,
we obtain
$
\alpha_2 = \left(\delta_1  + k\, \beta_1\,    \alpha_1\right) \alpha_1
< \left(\delta_0  + k\, \beta_0\,    \alpha_0\right) \alpha_0
= \alpha_1
$
and
$\tilde{\alpha_2}  =  \left({1+M}  + k\,  \alpha_{2} \right) \alpha_{2} <  \left({1+M}  + k\,  \alpha_{1} \right) \alpha_{1} =\tilde{\alpha_1}$.

To analyze the sequences $\{d_n\}$ and $\{\delta_n\}$, we must also have in mind that
$$
\beta^2_{1}  {\alpha_2} < (1+d_0)^2 \beta^2_0 \alpha_2
<
 \beta^2_0 (1+d_0)^2 \left(\delta_1  + k\, \beta_1\,    \alpha_1\right) \alpha_1
$$
$$ <  \beta^2_0 (1+d_0)^2 \left(\delta_0  + k\, \beta_0\,    \alpha_0\right) \alpha_1
 < \beta^2_0 \alpha_1
 $$
and therefore it follows that
 $\beta^2_{1}  \tilde{\alpha}_{2}  < \beta^2_0 \tilde{\alpha}_{1}$, $\beta_{1}  {\alpha_2} <  \beta_0 \alpha_1$
  and $\beta_{1}  \tilde{\alpha}_{2}  < \beta_0 \tilde{\alpha}_{1}$.

Then, to finish the case $n=2$, taking into account that by hypothesis $\delta_1<\delta_0$ , we get
$
 \delta_2= \delta^2_{1} +  k \, M \beta^2_{1}  (\alpha_{2} + \tilde{\alpha}_ {2})
   <\delta^2_0 +  k \, M \beta^2_0  (\alpha_{1} +  \tilde{\alpha}_{1}) =
 \delta_1
 $
 and
$
d_1=\delta_{1} + k \,  \beta_{1}  (\alpha_{2} +  \tilde{\alpha}_{2})
< \delta + k \,  \beta  (\alpha_{1} +  \tilde{\alpha}_{1}) = d_0.
$

From now, we suppose that
 $ \alpha_0>  \alpha_1 > \dots >  \alpha_n $, \,
 $\beta_0  {\alpha_0} > \beta_{1}  {\alpha_1}  >\dots >  \beta_{n}  \alpha_{n}$ and  $ \beta^2_0  {\alpha_1} > \beta^2_{1}  {\alpha_2}  >\dots >  \beta^2_{n}  \alpha_{n+1} $    hold,
 which implies that the sequences
 $\{\tilde{ \alpha}_k\}_{k=1}^{n}$, $\{\beta_{k}  \tilde{ \alpha}_{k}\}_{k=1}^{n}$,
 $\{\beta^2_k \, \tilde{ \alpha}_{k+1}\}_{k=1}^{n}$, $\{\beta_k \,  \alpha_{k+1}\}_{k=1}^{n}$
and  $\{\beta_k \,  \tilde{\alpha}_{k+1}\}_{k=1}^{n}$ are decreasing, as well as  $\{\delta_k \}_{k=1}^{n}$ and
  $\{d_k \}_{k=0}^{n-1}$ are.

 We need to check the inductive step.

 In first place, it is easy to prove:
 $$ \alpha_{n+1}= (\delta_n + k \,  \beta_{n} \alpha_{n}) \alpha_{n}
<  (\delta_{n-1} + k \,  \beta_{n-1} \alpha_{n-1}) \alpha_{n-1}=\alpha_n,
$$
$$
  \tilde{\alpha}_{n+1} = (1 + M +  k \,   \alpha_{n+1}) \alpha_{n+1}
  <  (1 + M +  k \,   \alpha_{n}) \alpha_{n} = \tilde{\alpha}_{n},
$$
$$
d_{n}=\delta_{n} + k \,  \beta_{n}  (\alpha_{n+1} +  \tilde{\alpha}_{n+1})
<
 \delta_{n-1} + k \,  \beta_{n-1}  (\alpha_{n} +  \tilde{\alpha}_{n})= d_{n-1},
$$
$$
\delta_{n+1}= \delta^2_{n} +  k \, M \beta^2_{n}  (\alpha_{n+1} +  \tilde{\alpha}_{n+1})
< \delta^2_{n-1} +  k \, M \beta^2_{n-1}  (\alpha_{n} +  \tilde{\alpha}_{n}) = \delta_n.
$$
On the other hand
$$
\beta_{n+1} \alpha_{n+1}  < \beta_{n}  \alpha_{n}
\Leftrightarrow
(1+d_{n}) \alpha_{n+1}  <  \alpha_{n}
\Leftrightarrow
(1+d_{n})   (\delta_{n} + k \,  \beta_{n} \alpha_{n}) \alpha_{n}  <  \alpha_{n}
$$
$$
\Leftrightarrow
(1+d_{n})   (\delta_{n} + k \,  \beta_{n} \alpha_{n})   < 1,
$$
which is true since that
$$
(1+d_{n})   (\delta_{n} + k \,  \beta_{n} \alpha_{n})   < (1+d_{0})   (\delta + k \,  \beta  \alpha_{0}) <1.
$$
Note that we also have
$$ \alpha_{n+2}= (\delta_{n+1} + k \,  \beta_{n+1} \alpha_{n+1}) \alpha_{n+1}
<  (\delta_{n} + k \,  \beta_{n} \alpha_{n}) \alpha_{n}=\alpha_{n+1}
$$
 and therefore,
$$
\beta^2_{n+1}  \alpha_{n+2}
 < (1+d_{n})^2 \beta^2_{n} \alpha_{n+2}
 <  \beta^2_{n} (1+d_{n})^2 \left(\delta_{n+1}  + k\, \beta_{n+1}\,    \alpha_{n+1}\right) \alpha_{n+1}
 $$
$$
  <  \beta^2_{n} (1+d_{n})^2 \left(\delta_{n}  + k\, \beta_{n}\,    \alpha_{n}\right) \alpha_{n}
 < \beta^2_{n} \alpha_{n+1}.
$$
Consequently $\{ \alpha_n\}$, $\{  \tilde{\alpha}_{n}\}$, $\{ d_n\}$ and $\{ \delta_n\}$
are  decreasing.

The proof is complete.
\qquad$\blacksquare$
\bigskip

%%%%%%%%%%%%%%%%%%%%%%%%%%%%%%%%%%%%%%%%%%%%%%
\subsection{A local convergence result}
%%%%%%%%%%%%%%%%%%%%%%%%%%%%%%%%%%%%%%%%%%%%%%

First of all, we notice that, fixed $x_0 \in B (x^*, r)$, if there is $\tilde{r}$, with the condition $(1+M+k r) r< \tilde{r}$ (in order to verify $B(x^*,\tilde{r})\subseteq \Omega$) and, in addition, $r$ verifies the two other conditions given in (\ref{cond}), then $\{x_n\}$ will be convergent.

On the other hand, the conditions given in (\ref{cond}) can be written in the following way:

\begin{equation} \label{cond1}
(1+M+k r) r < \tilde{r},
\end{equation}

\begin{equation} \label{cond2}
\delta_{1} < \delta_{0}
\Leftrightarrow
\delta^2+k M \beta^2 (\alpha_1+\tilde{\alpha}_{1})< \delta
\Leftrightarrow
0< \delta (1-\delta)- r P_1(r),
\end{equation}

$$
(1+d_0)^2 ( \delta_0 + k \beta \alpha_0 )<1
\Leftrightarrow
(1+\delta+k\beta(\alpha_1+\tilde{\alpha}_{1}))^2(\delta+k\beta r)<1
$$
\begin{equation} \label{cond3}
\Leftrightarrow
0< (1-(1+\delta)^2 \delta)-r P_3(r),
\end{equation}
where $P_1$ is a polynomial of degree one and $P_3$ is a polynomial of degree three, both decreasing and concave in $(0,+\infty)$. In this situation, it is clear that if $1-(1+\delta)^2 \delta>0$, then there is always $r$ checking (\ref{cond2}) and (\ref{cond3}).

 \begin{theo}\quad

Using the above notations,  under the initial conditions $(C1)-(C4)$, we assume that there exists $r>0$ verifying the conditions (\ref{cond1}),(\ref{cond2}), (\ref{cond3}) and $B(x^*,r)\subseteq \Omega$. Then, if we consider $x_0 \in B(x^*,r)$, the sequence $\{x_n\}$ given by (\ref{metUlm})  is well defined and converges to a solution $  x^{\ast}$  of $ F(x)=0$.
 \end{theo}
\textit{\textbf{Proof.}}\quad

First, it is easy to prove from the hypothesis that $x_n, x_n+F(x_n) \in  {\Omega}$ for $n\geq 1$. Then, the sequence $\{x_n\}$ given by (\ref{metUlm})  is well defined.

On the other hand, if we denote $L= \delta + k \beta r$, we have
$$
\Vert x_{1} - x^{*} \Vert
<  ( \delta + k \beta r )\Vert x_{0} - x{*} \Vert = L \Vert x_{0} - x^{*} \Vert
$$
$$
\Vert x_{2} - x^{*} \Vert
<  ( \delta_1 + k \beta_1  \alpha_1 )\Vert x_{1} - x{*} \Vert
<  ( \delta + k \beta r ) \Vert x_{1} - x^{*} \Vert < L^2 \Vert x_{0} - x^{*} \Vert
$$
So, in general, we obtain
$$
\Vert x_{n} - x^{*} \Vert
<  ( \delta_{n-1} + k \beta_{n-1}  \alpha_{n-1} )\Vert x_{n-1} - x{*} \Vert
<  ( \delta + k \beta r ) \Vert x_{n-1} - x^{*} \Vert < L^n \Vert x_{0} - x^{*} \Vert
$$
Since $L<1$ from lemma \ref{Lem2}, it follows that the sequence $\{x_n\}$ given by (\ref{metUlm}) converges to a solution $ x^{\ast}$  of $ F(x)=0$.
\qquad$\blacksquare$
\bigskip

Notice that, if $\Omega=\R^m$, the condition (\ref{cond1}) it is not necessary.

%%%%%%%%%%%%%% An example %%%%%%%%%%%%%%%%%%%%%%%%%
\subsection{An example}\label{exam}

Next, we illustrate the previous result with the following example given in~\cite{DS}.
We choose the max-norm.

Let $F:\mathbb{R}^3\rightarrow \mathbb{R}^3$ be defined as $F(x,y,z)=(x,y^2+y,e^z-1)$.
It is obvious that the unique solution of the system is $x^*=(0,0,0)$.

From $F$, having into account (\ref{DD}), we have
\[
F'(x,y,z)=\left(
\begin{array}{ccc}
1&0&0\\
0&2y+1&0\\
0& 0 &e^z
\end{array}
\right)
\]

\[
[(x,y,z),(u,v,w);F]=\left(
\begin{array}{ccc}
1&0&0\\
0&y+v+1&\\
0& 0 &{{e^z-e^w}\over{z-w} }
\end{array}
\right).
\]
So, $F'(x^*)$ is the identity matrix $3\times3$.
Then, $\| F'(x^*) \|=1$ and consider   $B_0= {\rm diag} \{\beta, \beta, \beta \}$  with $\beta = 0.75$ and $\delta = 0.25$.
On the other hand, there exists $\tilde{r}=1$, such that
$B(0,\tilde{r})=\{w\in\mathbb{R}^3:\,\|w\|<1\}\subset\mathbb{R}^3$,
and it is easy to prove that
$$
\|[(x,y,z),(u,v,w);F]-F'(0,0,0) \| \leq \max \{ \vert  y+v \vert,  \vert {{e^z-e^w}\over{z-w} } -1 \vert \}
$$
$$
 \leq   \|(x,y,z))\| +\|(u,v,w)\|,
$$
in $B(0,\tilde{r})$.
Therefore, $M=1, \, k =1$ and considering $\alpha_0= 0.246627$ we obtain:
$$
(1+M+k r) r = 0.554078< \tilde{r},
\quad
\alpha_1= 0.107275,
\quad
\tilde{\alpha_1} = 0.226058,
$$
$$
\delta - \delta_1 = 5.55112\times10^{-17},
\quad
d_0= 0.5,
\quad
1- (1+d_0)^2 ( \delta_0 + k \beta \alpha_0 ) = 0.0213177
$$
Therefore
the iterative process (\ref{metUlm}) is convergent from any starting point belonging to $B(x^*, 0.246627)$.

%%%%%%%%%%%%%%%%%%%%%%%%%%%%
\section{Numerical experiments}
%%%%%%%%%%%%%%%%%%%%%%%%%%%%%

In this section, we include two experiments to test the proposed algorithm  (\ref{metUlm}). In the first one, we check numerically its order of convergence and its stability behavior. Moreover, we propose specific chooses of the initial matrix $B_0$ in order to improve the stability of the classical Steffensen method (\ref{metStef}).
In the final example, we apply our iterative method to solve the nonlinear systems of equations that appear
to approximate a stiff differential problem with an implicit method.

\subsection{Academic example}

We consider the academic system of equations $F(x,y)=(0,0)$ given by:
\begin{eqnarray}\label{pro1}
(2x-\frac{x^2}{\epsilon})+(y-\frac{y^2}{2 \epsilon})&=&0,\\\nonumber
x+y&=&0.
\end{eqnarray}

This system has as solution $(x^*,y^*)=(0,0)$ and its Jacobian matrix is given by the following matrix
$$
\left(
  \begin{array}{cc}
    2-2\frac{x}{\epsilon} & 1-\frac{y}{\epsilon} \\
    1 & 1 \\
  \end{array}
\right)
$$

On the other hand, having into account (\ref{DD}), the matrix $[x,\tilde{x};F]$, with $x=(x_1,y_1)$ and $\tilde{x}=(x_2,y_2)$, is given by

$$
\left(
  \begin{array}{cc}
    2-\frac{x_1+x_2}{\epsilon} & 1-\frac{y_1+y_2}{\epsilon} \\
    1 & 1 \\
  \end{array}
\right)
$$

\begin{itemize}

\item Stability

The parameter $\epsilon$ (when $\epsilon \to 0$) increases the condition number of the difference divided matrix for a given initial guess  and the Steffensen method (\ref{metStef}) should have problems of convergence
to the solution for small values of $\epsilon$. Moreover, notice that $F^{'}(\epsilon,\epsilon)$ is not invertible.

We compute the maximum of the condition numbers of the linear systems that appear in the application of the Steffensen method (\ref{metStef}) ($||A||\cdot ||A^{-1}||$) and the maximum of the conditions numbers in all the matrix multiplications in the Moser-Steffensen type method (\ref{metUlm}) ($\frac{||A||\cdot ||B||}{||A B||}$).

%
%
%%%%%%%%%%%%%%%%%%%%%%%%%%%%%%%%%%%%%%%%%%%%%%%%%%%%%%%%%%%%%%%%%%%%%%%%%%%%%%%%%%%%%%%%%%%%%%%%%%%%%%%%%%%%%%%%%%%%
%%%%%%%%%%%%%%%%%%%%%%%%%%%%%%%%%%%%%%%%%%%%% begin Table %%%%%%%%%%%%%%%%%%%%%%%%%%%%%%%%%%%%%%%%%%%%%%%%%%%%%%%%%%
%
%
\begin{table}[ht!!!!]
\caption{System (\ref{pro1}). Errors for  $(x_0,y_0)=(-1,1)$,  $||B_0 F^{'}(-1,1)||\leq 10^{-3}$ and $\epsilon=1$.}
\vspace{-3mm}
\begin{center} {
\begin{tabular}{lccc}
\hline \hline
$n$ & Steffensen & Moser-Steffensen \\[0.01em]
\hline \\[-0.6em]
$1$ & $1.55 \ 10^{1}$  & $1.00 \ 10^{0}$ \\ \hline
$2$ & $1.84 \ 10^{1}$  & $3.76 \ 10^{-1}$ \\ \hline
$3$ & $2.12 \ 10^{1}$  & $1.08 \ 10^{-1}$ \\ \hline
$4$ & $2.41 \ 10^{1}$  & $1.80 \ 10^{-2}$ \\ \hline
$5$ & $2.69 \ 10^{1}$  & $9.12 \ 10^{-4}$ \\ \hline
$6$ & $2.98 \ 10^{1}$  & $3.61 \ 10^{-6}$ \\ \hline
$7$ & $3.26 \ 10^{1}$  & $7.60  \ 10^{-11}$ \\ \hline
$8$ & $3.54 \ 10^{1}$  & $4.18  \ 10^{-20}$ \\ \hline\hline
\end{tabular} }
\end{center}
\end{table}\label{table1}
%
%
%%%%%%%%%%%%%%%%%%%%%%%%%%%%%%%%%%%%%%%%%%%%%%%%%%%%%%%%%%%%%%%%%%%%%%%%%%%%%%%%%%%%%%%%%%%%%%%%%%%%%%%%%%%%%%%%%%%%
%%%%%%%%%%%%%%%%%%%%%%%%%%%%%%%%%%%%%%%%%%%%% end Table  %%%%%%%%%%%%%%%%%%%%%%%%%%%%%%%%%%%%%%%%%%%%%%%%%%%%%%%%%%%
%
%

%
%
%%%%%%%%%%%%%%%%%%%%%%%%%%%%%%%%%%%%%%%%%%%%%%%%%%%%%%%%%%%%%%%%%%%%%%%%%%%%%%%%%%%%%%%%%%%%%%%%%%%%%%%%%%%%%%%%%%%%
%%%%%%%%%%%%%%%%%%%%%%%%%%%%%%%%%%%%%%%%%%%%% begin Table %%%%%%%%%%%%%%%%%%%%%%%%%%%%%%%%%%%%%%%%%%%%%%%%%%%%%%%%%%
%
%
\begin{table}[ht!!!!]
\caption{System (\ref{pro1}). Errors for  $(x_0,y_0)=(-0.25,0.25)$,  $||B_0 F^{'}(-0.25,0.25)||\leq 10^{-3}$ and $\epsilon=10^{-1}$.}
\vspace{-3mm}
\begin{center} {
\begin{tabular}{lccc}
\hline \hline
$n$ & Steffensen & Moser-Steffensen \\[0.01em]
\hline \\[-0.6em]
$1$ & $1.70 \ 10^{0}$  & $6.84 \ 10^{-1}$ \\ \hline
$2$ & $1.98 \ 10^{0}$  & $2.08 \ 10^{-1}$ \\ \hline
$3$ & $2.27 \ 10^{0}$  & $8.90 \ 10^{-2}$ \\ \hline
$4$ & $2.55 \ 10^{0}$  & $3.95 \ 10^{-2}$ \\ \hline
$5$ & $2.83 \ 10^{0}$  & $6.80 \ 10^{-2}$ \\ \hline
$6$ & $3.12 \ 10^{0}$  & $6.95 \ 10^{-4}$ \\ \hline
$7$ & $3.41 \ 10^{0}$  & $1.22  \ 10^{-5}$ \\ \hline
$8$ & $3.69 \ 10^{0}$  & $5.38  \ 10^{-9}$ \\ \hline
$9$ & $3.97 \ 10^{0}$  & $1.33  \ 10^{-15}$ \\ \hline
$10$ & $4.25 \ 10^{0}$  & $1.00  \ 10^{-28}$ \\ \hline\hline
\end{tabular} }
\end{center}
\end{table}\label{table2}
%
%
%%%%%%%%%%%%%%%%%%%%%%%%%%%%%%%%%%%%%%%%%%%%%%%%%%%%%%%%%%%%%%%%%%%%%%%%%%%%%%%%%%%%%%%%%%%%%%%%%%%%%%%%%%%%%%%%%%%%
%%%%%%%%%%%%%%%%%%%%%%%%%%%%%%%%%%%%%%%%%%%%% end Table  %%%%%%%%%%%%%%%%%%%%%%%%%%%%%%%%%%%%%%%%%%%%%%%%%%%%%%%%%%%
%
%

In table 1, the  vector $(\epsilon,\epsilon)$ is inside of the ball containing the initial guess and the solution. The Steffensen method  (\ref{metStef}) diverges.

This numerical behavior is similar for smaller parameters of $\epsilon$, as we can see in the table 2 ($\epsilon=10^{-1}$). The maximum of the condition numbers for the  Steffensen method  (\ref{metStef}) is $6.09 \ 10^2$ (too big) and for the Moser-Steffensen type method (\ref{metUlm}) smaller than 30. Moreover, the sequences of condition numbers for both
 methods are increasing and decreasing sequences.

 For this example,
 the Steffensen method (\ref{metStef}) has a small region of convergence. However, the method (\ref{metUlm})
 only reduces a little its velocity.

This is the main advantage of the Moser-Steffensen type method (\ref{metUlm}). The condition number of the
operations used in the classical Steffensen method (\ref{metStef}) should be large (in the previous cases of divergence the condition number goes to infinity) and there is not a general strategy to find preconditioners for a given linear system. However, by construction, the condition number in the operations of the proposed method (\ref{metUlm}) (matrix multiplications) seems controlled with our election of $B_0$.

In Moser-type algorithms, as  (\ref{metUlm}), the sequence of matrices $B_n$  converges to the inverse of the Jacobian at the solution. For this reason, a good candidate for $B_0$ is an approximation to the inverse of the Jacobian at the initial guess \cite{GoLo}.

\item Order of convergence

In table 3, the method obtains the results expected by our theoretical analysis and the second order
convergence is clear. The vector $(\epsilon,\epsilon)$
is outside the convergence region associated to the initial guess. We compute the maximum of the condition numbers of the linear systems that appear in the application of the Steffensen method (\ref{metStef}) and the maximum of the conditions numbers in all the matrix multiplications in our method  (\ref{metUlm}). In this case, the maximum for the Steffensen method (\ref{metStef})  and for our method (\ref{metUlm}) are smaller than 10, and both methods work well.

%
%
%%%%%%%%%%%%%%%%%%%%%%%%%%%%%%%%%%%%%%%%%%%%%%%%%%%%%%%%%%%%%%%%%%%%%%%%%%%%%%%%%%%%%%%%%%%%%%%%%%%%%%%%%%%%%%%%%%%%
%%%%%%%%%%%%%%%%%%%%%%%%%%%%%%%%%%%%%%%%%%%%% begin Table %%%%%%%%%%%%%%%%%%%%%%%%%%%%%%%%%%%%%%%%%%%%%%%%%%%%%%%%%%
%
%
\begin{table}[ht!!!!]\label{ta1}
\caption{System (\ref{pro1}). Errors for $(x_0,y_0)=(-1,1)$,  $||B_0 F^{'}(-1,1)||\leq 10^{-3}$ and $\epsilon=3$.}
\vspace{-3mm}
\begin{center} {
\begin{tabular}{lccc}
\hline \hline
$n$ & Steffensen & Moser-Steffensen \\[0.01em]
\hline \\[-0.6em]
$1$ & $1.41 \ 10^{0}$  & $3.55 \ 10^{-1}$ \\ \hline
$2$ & $7.61 \ 10^{-1}$  & $5.09 \ 10^{-2}$ \\ \hline
$3$ & $2.40 \ 10^{-1}$  & $1.86 \ 10^{-3}$ \\ \hline
$4$ & $2.60 \ 10^{-2}$  & $3.74 \ 10^{-6}$ \\ \hline
$5$ & $3.11 \ 10^{-4}$  & $2.01 \ 10^{-11}$ \\ \hline
$6$ & $4.56 \ 10^{-8}$  & $7.10 \ 10^{-22}$ \\ \hline
$7$ & $9.79 \ 10^{-16}$  & \\ \hline \hline
\end{tabular} }
\end{center}
\end{table}\label{table3}
%
%
%%%%%%%%%%%%%%%%%%%%%%%%%%%%%%%%%%%%%%%%%%%%%%%%%%%%%%%%%%%%%%%%%%%%%%%%%%%%%%%%%%%%%%%%%%%%%%%%%%%%%%%%%%%%%%%%%%%%
%%%%%%%%%%%%%%%%%%%%%%%%%%%%%%%%%%%%%%%%%%%%% end Table  %%%%%%%%%%%%%%%%%%%%%%%%%%%%%%%%%%%%%%%%%%%%%%%%%%%%%%%%%%%
%
%

%
%
%%%%%%%%%%%%%%%%%%%%%%%%%%%%%%%%%%%%%%%%%%%%%%%%%%%%%%%%%%%%%%%%%%%%%%%%%%%%%%%%%%%%%%%%%%%%%%%%%%%%%%%%%%%%%%%%%%%%
%%%%%%%%%%%%%%%%%%%%%%%%%%%%%%%%%%%%%%%%%%%%% begin Table %%%%%%%%%%%%%%%%%%%%%%%%%%%%%%%%%%%%%%%%%%%%%%%%%%%%%%%%%%
%
%
\begin{table}[ht!!!!]
\caption{System (\ref{pro1}). Errors for $(x_0,y_0)=(-0.5,0.5)$,  $||B_0 F^{'}(-0.5,0.5)||\leq 10^{-3}$ and $\epsilon=1$.}
\vspace{-3mm}
\begin{center} {
\begin{tabular}{lccc}
\hline \hline
$n$ & Steffensen & Moser-Steffensen \\[0.01em]
\hline \\[-0.6em]
$1$ & $1.49 \ 10^{0}$  & $2.12 \ 10^{-1}$ \\ \hline
$2$ & $2.75 \ 10^{0}$  & $4.21 \ 10^{-2}$ \\ \hline
$3$ & $3.35 \ 10^{0}$  & $3.09 \ 10^{-3}$ \\ \hline
$4$ & $6.82 \ 10^{0}$  & $2.68 \ 10^{-5}$ \\ \hline
$5$ & $9.80 \ 10^{0}$  & $2.77 \ 10^{-9}$ \\ \hline
$6$ & $1.27 \ 10^{1}$  & $3.75 \ 10^{-17}$ \\ \hline
$7$ & $1.56 \ 10^{1}$  & $8.71  \ 10^{-33}$ \\ \hline
$8$ & $1.84 \ 10^{1}$  & \\ \hline
$9$ & $2.13 \ 10^{1}$  &  \\ \hline
$10$ & $2.41 \ 10^{1}$  & \\ \hline
 \hline
\end{tabular} }
\end{center}
\end{table}\label{table4}
%
%
%%%%%%%%%%%%%%%%%%%%%%%%%%%%%%%%%%%%%%%%%%%%%%%%%%%%%%%%%%%%%%%%%%%%%%%%%%%%%%%%%%%%%%%%%%%%%%%%%%%%%%%%%%%%%%%%%%%%
%%%%%%%%%%%%%%%%%%%%%%%%%%%%%%%%%%%%%%%%%%%%% end Table  %%%%%%%%%%%%%%%%%%%%%%%%%%%%%%%%%%%%%%%%%%%%%%%%%%%%%%%%%%%
%
%

The results in the table 4 are similar to the first case, both methods have second order of convergence.

In table 4 case, we consider a smaller $\epsilon$ but again the point $(\epsilon,\epsilon)$ is outside the ball including the initial guess (closer also to the solution) and the solution. In this case, our method
preserve the second order of convergence (condition number smaller than 10), however the Steffensen method diverges (with condition number $1.80 \ 10^2$ in the last iteration).

\item Election of the initial matrix $B_0$

 In general, as we indicated before, a good candidate for $B_0$ is an approximation to the inverse of the Jacobian matrix at the initial guess.

On the other hand,  it is not necessary to consider  $B_0$
 a really accurate approximation to the Jacobian at the initial guess, as we can see in table 5. In particular, we can take some iterations of some of the
 algorithms that appear for instance in \cite{AEH}.

%
%
%%%%%%%%%%%%%%%%%%%%%%%%%%%%%%%%%%%%%%%%%%%%%%%%%%%%%%%%%%%%%%%%%%%%%%%%%%%%%%%%%%%%%%%%%%%%%%%%%%%%%%%%%%%%%%%%%%%%
%%%%%%%%%%%%%%%%%%%%%%%%%%%%%%%%%%%%%%%%%%%%% begin Table %%%%%%%%%%%%%%%%%%%%%%%%%%%%%%%%%%%%%%%%%%%%%%%%%%%%%%%%%%
%
%
\begin{table}[ht!!!!]
\caption{System (\ref{pro1}). Errors for Moser-Steffensen, $(x_0,y_0)=(-2,2)$ and $\epsilon=3$. }
\vspace{-3mm}
\begin{center} {
\begin{tabular}{lccc}
\hline \hline
$n$ & $||B_0 F^{'}(-2,2)||\leq 1$ & $||B_0 F^{'}(-2,2)||\leq 10^{-1}$ & $||B_0 F^{'}(-2,2)||\leq 10^{-3}$ \\[0.01em]
\hline \\[-0.6em]
$1$ & $1.84 \ 10^{0}$ & $1.10 \ 10^{0}$ & $9.44 \ 10^{-1}$   \\ \hline
$2$ & $8.75  \ 10^{-1}$ & $2.83 \ 10^{-1}$ & $2.24 \ 10^{-1}$   \\ \hline
$3$ & $2.53  \ 10^{-1}$ & $3.56 \ 10^{-2}$ & $2.43 \ 10^{-2}$   \\ \hline
$4$ & $3.29  \ 10^{-2}$ & $9.86 \ 10^{-4}$ & $4.87 \ 10^{-4}$  \\ \hline
$5$ & $8.61 \ 10^{-4}$ & $1.12 \ 10^{-6}$ & $2.82 \ 10^{-7}$   \\ \hline
$6$ & $8.37  \ 10^{-7}$ & $1.89 \ 10^{-12}$ & $1.21 \ 10^{-13}$   \\ \hline
$7$ & $1.01  \ 10^{-12}$ & $6.57 \ 10^{-24}$ & $2.71 \ 10^{-26}$   \\ \hline
$7$ & $1.83  \ 10^{-24}$ &  &    \\ \hline\hline
\end{tabular} }
\end{center}
\end{table}\label{ta5}
%
%
%%%%%%%%%%%%%%%%%%%%%%%%%%%%%%%%%%%%%%%%%%%%%%%%%%%%%%%%%%%%%%%%%%%%%%%%%%%%%%%%%%%%%%%%%%%%%%%%%%%%%%%%%%%%%%%%%%%%
%%%%%%%%%%%%%%%%%%%%%%%%%%%%%%%%%%%%%%%%%%%%% end Table  %%%%%%%%%%%%%%%%%%%%%%%%%%%%%%%%%%%%%%%%%%%%%%%%%%%%%%%%%%%
%
%

Finally, in the table 7,
we force the method to take the bad iteration $(\epsilon, \epsilon)$ (here the Jacobian is not invertible). Only with the new approach we are able to find
the solution (the Steffensen method diverges).

Indeed, we consider $(x_0,y_0)=(\epsilon,\epsilon)$. For this situations, a possibility is to take the initial matrix
$B_0=\delta I_2$, where $\delta$ is a small parameter ($10^{-2}$
in our experiments) and
$I_2$ is the identity matrix. When the method leaves the conflict zone, it recovers its good properties
(second order convergence).

%
%
%%%%%%%%%%%%%%%%%%%%%%%%%%%%%%%%%%%%%%%%%%%%%%%%%%%%%%%%%%%%%%%%%%%%%%%%%%%%%%%%%%%%%%%%%%%%%%%%%%%%%%%%%%%%%%%%%%%%
%%%%%%%%%%%%%%%%%%%%%%%%%%%%%%%%%%%%%%%%%%%%% begin Table %%%%%%%%%%%%%%%%%%%%%%%%%%%%%%%%%%%%%%%%%%%%%%%%%%%%%%%%%%
%
%
\begin{table}[ht!!!!]
\caption{System (\ref{pro1}). Errors for $(x_0,y_0)=(2,2)$, $\epsilon=2$ and $B_0= 10^{-2} \ I_2$. }
\vspace{-3mm}
\begin{center} {
\begin{tabular}{lcc}
\hline \hline
$n$ & Moser-Steffensen \\[0.01em]
\hline \\[-0.6em]
$10$ & $1.13 \ 10^{-2}$ \\ \hline
$11$ & $2.81 \ 10^{-4}$ \\ \hline
$12$ & $2.07  \ 10^{-7}$ \\ \hline
$13$ & $1.30  \ 10^{-13}$ \\ \hline
$14$ & $5.88  \ 10^{-26}$ \\ \hline \hline
\end{tabular} }
\end{center}
\end{table}\label{ta7}
%
%
%%%%%%%%%%%%%%%%%%%%%%%%%%%%%%%%%%%%%%%%%%%%%%%%%%%%%%%%%%%%%%%%%%%%%%%%%%%%%%%%%%%%%%%%%%%%%%%%%%%%%%%%%%%%%%%%%%%%
%%%%%%%%%%%%%%%%%%%%%%%%%%%%%%%%%%%%%%%%%%%%% end Table  %%%%%%%%%%%%%%%%%%%%%%%%%%%%%%%%%%%%%%%%%%%%%%%%%%%%%%%%%%%
%
%

\end{itemize}

%%%%%%%%%%%%%%%%%%%%%%%%%%%%
%================
\subsection{A stiff problem: Chapman atmosphere}
%================
%%%%%%%%%%%%%%%%%%%%%%%%%%%%

This model represents the Chapman mechanism for the generation of the
ozone and the oxygen singlet. In this example, the concentration
of the oxygen $y_3=[O_2]$  will be held constant. It is a severe test
for a stiff ODE package \cite{Ha} governed by the following
equations:
\begin{eqnarray*}
y^{\prime }_{1}(t)&=&2 k_3(t) y_{3} + k_4(t) y_2(t) - (  k_1 y_3 + k_2 y_2(t)) y_{1}(t) ,\\
y^{\prime }_{2}(t)&=&k_1 y_1(t) y_3-(k_2 y_1(t)+k_4(t)) y_2(t), \\
\end{eqnarray*}
with $y_3=3.7 \times10^{16}$, $k_1=1.63 \times10^{-16}$, $k_2=4.66 \times
10^{-16}$,
\begin{equation*}
k_i(t)=\left\lbrace
\begin{array}{ll}
\exp (\frac{a_i}{\sin (\omega t)}), &\text{ if } \sin (\omega t)>0 \\
 0, & \text{ otherwise} \\
\end{array}
\right.
\end{equation*}
for $i=3,4,$ with $a_3=22.62$, $a_4=7.601$ and
$\omega=\frac{\pi}{43200}$. The constant $43200$ is $12$ h
measured in seconds.
 The initial conditions are $y_{1}(0)=10^6$ and $y_{2}(0)=10^{12}$.

This problem has important features like:
\begin{itemize}
\item
 The Jacobian matrix is not a
constant. \item The diurnal effect is present. \item The
oscillations are fast. \item The time interval used is fairly
long, $0\leq t \leq 8.64 \ 10^5$, or $10$ days.
\end{itemize}

Let $h>0$. Given different coefficients $c_i$, $1 \leq i \leq s$
there is a (unique for $h$ sufficiently small) polynomial of
collocation $q(t)$ of degree less than or equal to  $s$ such that
\begin{equation}\label{colocacion}
q(t_0)=y_0, \quad q'(t_0+c_i \, h)=f(t_0+c_i \, h,q(t_0+c_i \, h))
\quad \textnormal{ if }1 \leq i \leq s.
\end{equation}
The collocation methods are defined  by an approximation $y(t)
\simeq q(t)$, and are equivalent to implicit RK methods
of $s$ stages
\begin{equation}
\begin{array}{rcl}
k_i&=&f(t_0+c_i \, h,y_0+h\displaystyle\sum_{j=1}^s a_{i,j} \,
k_j), \\ [0.25cm] y_1&=&y_0+h \displaystyle\sum_{i=1}^s b_i \,
k_i,
\end{array}
\end{equation}

 \begin{figure}
\centering
\includegraphics[width=7cm,height=7cm]{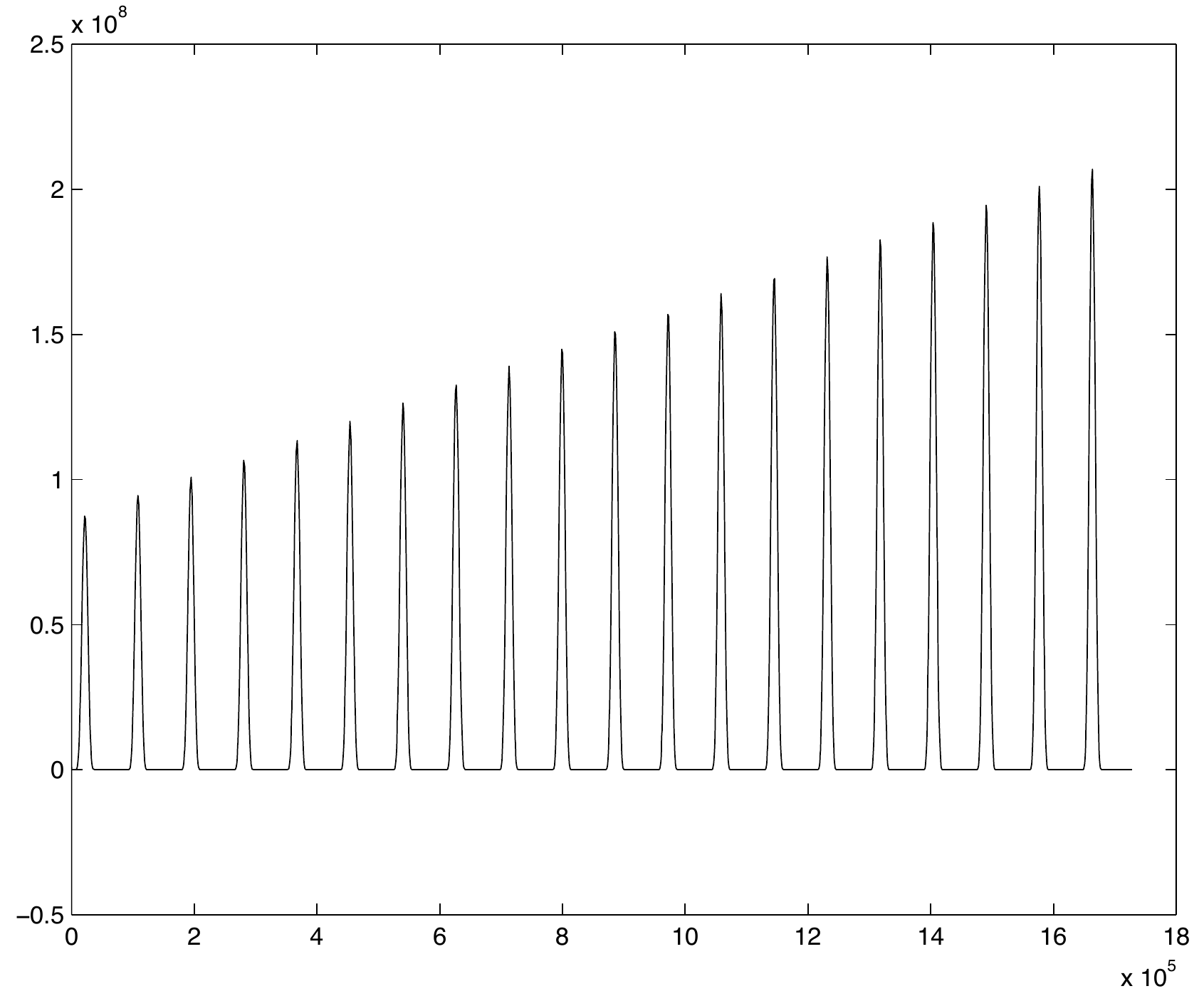}
\caption{First component of the Chapman atmosphere problem} \label{chapman1}
\end{figure}

for the coefficients
\begin{equation}\label{coeficientes-RK}
\begin{array}{rcl}
a_{i,j}&=&\displaystyle\int_0^{c_i} \displaystyle\prod_{l \neq
j}\frac{u-c_l}{c_j-c_l} \, du, \\ [0.5cm]
b_i&=&\displaystyle\int_0^1 \displaystyle\prod_{l \neq
i}\frac{u-c_l}{c_i-c_l} \, du.
\end{array}
\end{equation}
The coefficients $c_i$ play the role of the nodes of the
quadrature formula, and the associated coefficients $b_i$ are
analogous to the weights. From (\ref{coeficientes-RK}) we
can find implicit RK methods called Gauss of order $2s$,
Radau IA and Radau IIA of order $2s-1$ and Lobatto IIIA of order
$2s-2$. See \cite{Ha} for more details.

 \begin{figure}
\centering
\includegraphics[width=7cm,height=7cm]{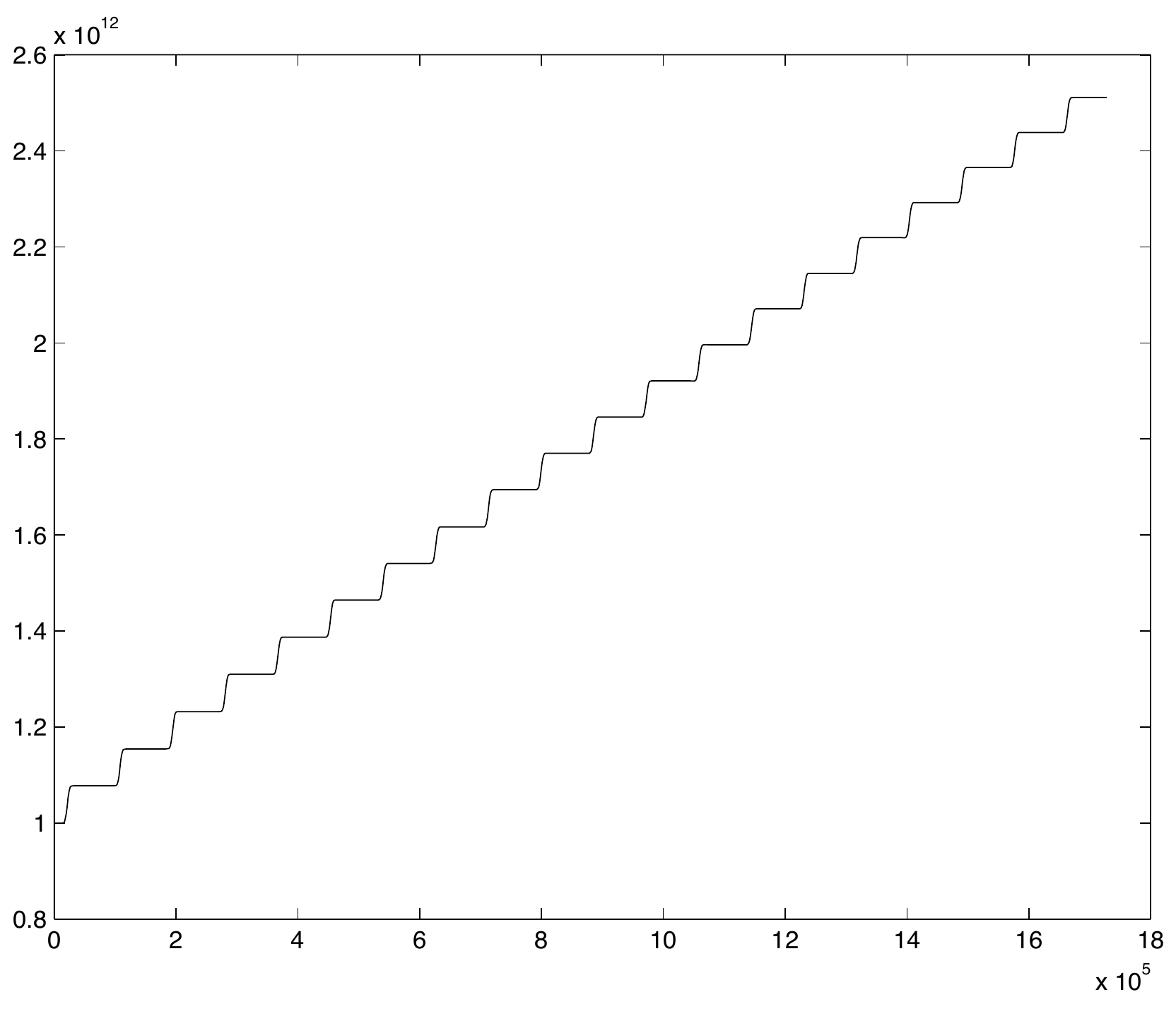}
\caption{Second component of the Chapman atmosphere problem} \label{chapman2}
\end{figure}

We consider the implicit fourth order Gauss
method ($s=2$ as collocation method). We approximate the associated nonlinear systems
of equations using our Moser-Steffensen's method as a black box.
We obtain a good approximation as we can see in Figures \ref{chapman1} and \ref{chapman2}. Note
that $y_2 = [0_3]$ looks like a staircase with a rise at midday
every day and $y_1 = [O]$ looks like a spike with its amplitude
increases each day.

%\clearpage

%%%%%%%%%%%%%%%%%%%%%%%%%%%%%%%%%%%%%%%%%%%%%%
%%%%%%%%%%%%%%%%%%%%%%%%%%%%%%%%%%%%%%%%%%%%%%

\end{document}